\documentclass[10pt]{article}
\usepackage{amssymb,latexsym,amsmath,epsfig,amsthm} %% Add other packages as necessary
\usepackage{mathtools}
\usepackage{array}

\makeatletter

\renewcommand\section{\@startsection {section}{1}{\z@}
{-30pt \@plus -1ex \@minus -.2ex}
{2.3ex \@plus.2ex}
{\normalfont\normalsize\bfseries}}

\renewcommand\subsection{\@startsection{subsection}{2}{\z@}
{-3.25ex\@plus -1ex \@minus -.2ex}
{1.5ex \@plus .2ex}
{\normalfont\normalsize\bfseries}}

\renewcommand{\@seccntformat}[1]{\csname the#1\endcsname. }

\makeatother

\numberwithin{equation}{section}

\newtheorem{thm}{Theorem}
\newtheorem{lem}{Lemma}

\numberwithin{thm}{section}
\numberwithin{prop}{section}
\numberwithin{col}{thm}
\numberwithin{lem}{section}

\begin{document}

\begin{center}
\uppercase{\bf Bernoulli-Stirling numbers}
\vskip 20pt
{\bf Ren\'{e} Gy}\\
{\tt rene.gy@numericable.com}\\
\end{center}
\vskip 20pt

%\centerline{\smallit Received: , Revised: , Accepted: , Published: } % We will fill in the dates
\vskip 30pt

\centerline{\bf Abstract}

\noindent Congruences modulo prime powers involving generalized Harmonic numbers are known. While looking for similar congruences, we have encountered a curious triangular array of numbers indexed with positive integers $n,k$, involving the Bernoulli and cycle Stirling numbers. These numbers are all integers and they vanish when $n-k$ is odd. This triangle has many similarities with the Stirling triangle. In particular, we show how it can be extended to negative indices and how this extension produces a {\it second kind} of such integers which may be considered as a new generalization of the Genocchi numbers and for which a generating function is easily obtained. But our knowledge of these integers remains limited, especially for those of the {\it first kind}. 

\pagestyle{myheadings} 
%%\markright{\smalltt INTEGERS: 17 (2017)\hfill} 
\thispagestyle{empty} 
\baselineskip=12.875pt 
\vskip 30pt

%11111111111111111111111111111111111111111111111111111111111111111111111111111111111111111111111111111111111111111111111111111111111111111111111111111111111111111111111111111111111111111111111111111111111111%
\section{Introduction}
Let $n$ and $k$ be non-negative integers and let the generalized Harmonic numbers $H_n^{(k)}$, $G_n^{(k)}$ be defined as
\begin{align*} H_n^{(k)}&:=\sum_{j=1}^n \frac{1}{j^k}, \\ G_n^{(k)}&:= \sum_{1\le i_1<i_2<\cdot\cdot <i_k\le n} \frac{1}{i_1 i_2 \cdot \cdot \cdot i_k}\text{,} \end{align*}
with $H_n^{(1)}=G_n^{(1)}=\sum_{j=1}^n \frac{1}{j} =H_n=G_n$; $H_n^{(0)}=n$ and $G_n^{(0)}=1$. We have \cite{Katriel10} $${n+1\brack k+1}= n! G_n^{(k)},$$ ${n\brack {k}}$ being the cycle Stirling number (or unsigned Stirling number of the first kind), so that the Harmonic and cycle Stirling numbers are inter-related by the convolution 
\begin{align} \label{conv} k{n+1\brack k+1}&= -\sum_{j=0}^{k-1}(-1)^{k-j}H_n^{(k-j)}{n+1\brack j+1}\text{,} \end{align}
which is obtained as a direct application of the well-known relation between elementary symmetric polynomials and power sums \cite{Kalman00}.
\noindent Extended congruences for the Harmonic numbers $H_{p-1}^{(k)}$, modulo any power of a prime $p$ are known \cite{Washington98}, \cite{Gy19}. Our initial motivation for the work reported in the present paper is to look for similar congruences modulo prime powers, involving $G_{p-1}^{(k)}$, or the cycle Stirling numbers ${p\brack k+1}$, instead of $H_{p-1}^{(k)}$. We will show that such similar congruences for $G_{p-1}^{(k)}$ do exist, but that they are just the particular prime instances of not very well-known but elementary identities for the cycle Stirling numbers. This will lead us to introduce a triangular array of integers, the {\it Bernoulli-Stirling numbers}, involving the Bernoulli and Stirling numbers which we believe is new.

%22222222222222222222222222222222222222222222222222222222222222222222222222222222222222222222222222222222222222222222222222222222222222222222222222222222222222222222222222222222222222222222222222222222222222%
\section{Notation and preliminaries}
In addition to what was exposed in the previous introduction, further notation that we use throughout this paper is presented in this section, along with classical results which we will need. Most of these results can be found in textbooks like \cite{Graham94} and they are given hereafter without proof.\\
\indent In the following, $p$ denotes a prime number, $x$ denotes the argument of a generating function, of a polynomial function or of a formal power series, $[[x^n]]f(x)$ the coefficient of $\frac{x^n}{n!}$ in $f(x)$ and $D^mf(x)$ is the $m$-order derivative of $f(x)$ with respect to $x$. We will use the Iverson bracket notation: $\big[\mathfrak{P}\big]=1$ when proposition $\mathfrak{P}$ is true, and  $\big[\mathfrak{P}\big]=0$ otherwise.
For $q$ a rational number, we denote by $v_p(q)$ the $p$-adic order of $q$ and we have $v_p(q_1\cdot q_2)=v_p(q_1)+v_p(q_2)$ and $v_p(q_1+ q_2)\ge \min{(v_p(q_1),v_p(q_2))}$.\\
The  binomial coefficients  ${n\choose k}$ are defined by  $\sum_{k}{n\choose k}x^k =(1+x)^n$,  whatever the sign of integer $n$. They obviously vanish when $k<0$. They are easily obtained by  the basic recurrence relation ${{n}\choose {k}}={{n-1}\choose {k}}+{{n-1}\choose{k-1}}$. When $n  > 0$, we have ${-n\choose k}= (-1)^k{n+k-1\choose n-1}$ and we have the well-known inversion formula
\begin{align}\label{hs1} \sum_{k\ge 0} (-1)^{k-j}{k\choose j}{n\choose k}&= [n=j].\end{align}
\indent The  cycle Stirling numbers  ${n\brack k}$, $n\ge0$ may be defined by the horizontal generating function
\begin{align}\label{gf}  \sum_{k\ge 0}{n\brack {k}}x^k&=\prod_{j=0}^{n-1}(x+j), \end{align}
where an empty product is meant to be $1$. Alternatively, they have the exponential generating function 
\begin{align} \label{egfs1}
 \sum_{n\ge 0} {n \brack k} \frac{x^n}{n!}&=\frac {(-1)^k \big(\ln(1-x)\big)^k}{k!}.\end{align}
They obviously vanish when $k<0$ and $k>n$. They are easily obtained by the basic recurrence ${{n}\brack {k}}=(n-1){{n-1}\brack {k}}+{{n-1}\brack {k-1}}$, valid for $n \ge 1$, with ${{0}\brack {k}}=[k=0]$. They also obey the generalized recurrence relation 
\begin{align} \label{recst1} {n+1\brack m+1}&= \sum_{h\ge 0}\binom{h+m}{m}{n\brack h+m}. \end{align}
We let  ${n\brace {k}}$, $n\ge0$, be the partition Stirling numbers (or Stirling numbers of the second kind). They also vanish when $k<0$ and $k>n$. Their basic recurrence is
${{n}\brace {k}}=k{{n-1}\brace {k}}+{{n-1}\brace {k-1}}$ for $n\ge1$, with $ {0\brace {k}}=[k=0]$. 
\noindent They have the following exponential generating function
\begin{align} \label{egfs2}  \sum_{n\ge 0}  {n \brace k} \frac{x^n}{n!}&=\frac {(e^x-1)^k}{k!}.\end{align}
We will make use of the following Lemma about the Stirling numbers, for which there is a proof in \cite{Graham94}, p. 266-271. 
 \begin{lem}\label{lemn-k}  Let $n, k$ be non-negative integers. There exists a polynomial  $Q_k \in {\mathbb Q}[X]$ of degree $2k$, such that $Q_k(n)$ coincides with ${n \brack n-k}$. If $k>0$ then $0,1,\cdot \cdot \cdot, k$ are roots of the polynomial $Q_k$. Moreover, we have $Q_k(-n)={n+k \brace n}$ and in particular $Q_k(-1)=1$. \end {lem}

\indent For $n\ge0$, let $B_n$ be the  Bernoulli numbers. The first of them are $B_0=1, B_1=-\frac{1}{2}, B_2=\frac{1}{6}, B_3=0, B_4=-\frac{1}{30},...$  and for $h>0$, we have   $B_{2h+1}=0$. They have the well-known exponential generating function
\begin{align}
\label{egfbern} \frac{t}{e^t-1}&=\sum_{n\ge0}B_n\frac{t^n}{n!} \end{align}
and they obey the recurrence
\begin{align}
\label{e2v} (-1)^nB_n&=\sum_{k=0}^{n}{n \choose k}B_k.  \end{align}
We will also make use of the Von Staudt-Clausen theorem which states that the denominator of $B_n$ in reduced form, is the product of all primes $p$ such that $p-1$ divides $n$. In particular, any prime may divide the denominator of a Bernoulli number once at most.

%333333333333333333333333333333333333333333333333333333333333333333333333333333333333333333333333333333333333333333333333333333333333333333333333333333333333333333333333333333333333333333333%
\section {Two identities for the cycle Stirling numbers}
\noindent In this section, we re-demonstrate two identities for the cycle Stirling numbers which, in spite of their similarity to Equation \eqref{recst1}, do not seem to be very well-known.
\begin {thm}\label{recst}
Let $m,n$ be non-negative integers. We have
\begin{equation} \label{recst2} {n+1\brack m+1}= (-1)^{n-m}\sum_{h\ge 0}\binom{h+m}{m}{n\brack h+m}(-n)^h \text{.}\end{equation}
Moreover, if $n>0$, 
\begin{equation} \label{recst3} {n\brack m}= (-1)^{n-m}\sum_{h\ge 0}\binom{h+m-1}{m-1}{n\brack h+m}(-n)^h \text{.} \end{equation}
In more compact and symmetric formulations, these two identies also read
\begin{align} \label{} (-1)^{n-m}{n+1\brack m+1}(-n)^m&=\sum_{h}\binom{h}{m}{n\brack h}(-n)^h \end{align}
and, if $n>0$,
\begin{align} \label{recst4} (-1)^{n-m}{n\brack m}(-n)^{m}&= \sum_{h}\binom{h-1}{m-1}{n\brack h}(-n)^{h} \text{.}\end{align}
\end{thm}
\noindent {\bf Remark.} These identities are not in \cite{Graham94} where quite many finite sums, recurrences and convolutions involving Stirling numbers are reported. Our Equation \eqref{recst2} may be obtained as a particular case of Theorem 3 in \cite{Agoh10}. An identity equivalent to our Equation \eqref{recst3}, is obtained incidentally in \cite{Adamchik97}, where it is not even labelled. Another identity, equivalent to our Equation \eqref{recst3} is the equation (18) from \cite{Shevelev12}, where it is said to be new.\\
\ \\
	{\it Proof of Theorem 3.1.} Like in \cite{Adamchik97}, our proof will highlight that Equation \eqref{recst2} and  Equation \eqref{recst3} are actually closely related to the convolution Equation \eqref{conv} between the Harmonic and cycle Stirling numbers. 
	Let $f_n(x):=\prod_{h=0}^{n-1}(x-h)$. We are going to show that, for $m\ge 1$,
	\begin{align} \label{conv2} m\frac{D^mf_n(x)}{m!}&=-\sum_{h=0}^{m-1}(-1)^{m-h}\frac{D^hf_n(x)}{h!}\sum_{j=0}^{n-1}\frac{1}{(x-j)^{m-h}} \text{.}\end{align}
	It is true for $m=1$, since $Df_n(x)=\sum_{j=0}^{n-1}\prod_{h\neq j}(x-h)=f_n(x)\sum_{j=0}^{n-1}\frac{1}{x-j}$. We suppose (induction hypothesis) that it is true for some $m$, so that 
	\begin{align*} (m+1)\frac{D^{m+1}f_n(x)}{(m+1)!}&=D\frac{D^mf_n(x)}{m!}=\frac{1}{m}D\left(m\frac{D^mf_n(x)}{m!}\right)\end{align*}
	\begin{align*}(m+1)\frac{D^{m+1}f_n(x)}{(m+1)!}&=\frac{1}{m}\left( -\sum_{h=0}^{m-1}(-1)^{m-h}\frac{D^{h+1}f_n(x)}{h!}\sum_{j=0}^{n-1}\frac{1}{(x-j)^{m-h}}\right)\\
	&+\frac{1}{m}\left(\sum_{h=0}^{m-1}(-1)^{m-h}\frac{D^hf_n(x)}{h!}(m-h)\sum_{j=0}^{n-1}\frac{1}{(x-j)^{m+1-h}}\right)\\
	&=\frac{1}{m}\left(\sum_{h=1}^{m}(-1)^{m-h}\frac{D^{h}f_n(x)}{(h-1)!}\sum_{j=0}^{n-1}\frac{1}{(x-j)^{m+1-h}}\right)\\
	&+\frac{1}{m}\left(m\sum_{h=0}^{m-1}(-1)^{m-h}\frac{D^hf_n(x)}{h!}\sum_{j=0}^{n-1}\frac{1}{(x-j)^{m+1-h}}\right)\\
	&-\frac{1}{m}\left(\sum_{h=1}^{m-1}(-1)^{m-h}\frac{D^hf_n(x)}{(h-1)!}\sum_{j=0}^{n-1}\frac{1}{(x-j)^{m+1-h}}\right)\\
	&=\frac{1}{m}\left(\frac{D^{m}f_n(x)}{(m-1)!}\sum_{j=0}^{n-1}\frac{1}{x-j}\right)\\
	&+\frac{1}{m}\left(m\sum_{h=0}^{m-1}(-1)^{m-h}\frac{D^hf_n(x)}{h!}\sum_{j=0}^{n-1}\frac{1}{(x-j)^{m+1-h}}\right)\\
	&=\sum_{h=0}^{m}(-1)^{m-h}\frac{D^hf_n(x)}{h!}\sum_{j=0}^{n-1}\frac{1}{(x-j)^{m+1-h}} \text { .}
	\end{align*}
	This establishes the validity of Equation \eqref{conv2}. Now, when $x=n$, it reads
	\begin{align} \label{conv3} m\frac{D^mf_n(n)}{m!}&=-\sum_{h=0}^{m-1}(-1)^{m-h}\frac{D^hf_n(n)}{h!}H_n^{(m-h)} \text{.}\end{align}
	This is the same recurrence as in the convolution Equation \eqref{conv}, with the same inital value, since by definition $f_n(n)=n!={n+1\brack 1}$. Then 
	$$\frac{D^mf_n(n)}{m!}= {n+1\brack m+1} \text{.}$$
	On the other hand, from Equation \eqref{gf}, we have
	$f_n(x)=\sum_{h=0}^n{n \brack h}(-1)^{n-h}x^h$, and hence
	\begin{align*} D^mf_n(x)&=m!\sum_{h=0}^n{n \brack h}{h \choose m}(-1)^{n-h}x^{h-m}\text {,}\\
	\frac{D^mf_n(n)}{m!}&=\sum_{h=0}^n{n \brack h}{h \choose m}(-1)^{n-h}n^{h-m}\text {.}\end{align*}
	Hence 
	\begin{align*} {n+1\brack m+1}&=\sum_{h=0}^n{n \brack h}{h \choose m}(-1)^{n-h}n^{h-m}\\
	&=\sum_{h=0}^{n-m}(-1)^{n-m}{n \brack h+m}{h+m \choose m}(-n)^{h} \text{.}
	\end{align*}
	This completes the proof of Equation \eqref{recst2}. For the proof of Equation \eqref{recst3}, we also use an induction argument, but on $m$ and backward. 
	Our induction hypothesis is \begin{align*} {n\brack m+1}&= (-1)^{n-(m+1)}\sum_{h=0}^{n-(m+1)}\binom{h+m}{m}{n\brack h+m+1}(-n)^h \text{.}\end{align*}
	Then \begin{align*} {n\brack m+1}&=- (-1)^{n-m}\sum_{h=1}^{n-m}\binom{h+m-1}{m}{n\brack h+m}(-n)^{h-1} \text{.}\end {align*}
	Hence \begin{align*} n{n\brack m+1}&= (-1)^{n-m}\sum_{h=1}^{n-m}\binom{h+m-1}{m}{n\brack h+m}(-n)^{h} \text{.} \end {align*}
	We subtract the latter equation from Equation \eqref{recst2}, and we obtain
	\begin{align*} {n+1\brack m+1} -n{n\brack m+1}&= (-1)^{n-m}\sum_{h=1}^{n-m}\left(\binom{h+m}{m}-\binom{h+m-1}{m}\right){n\brack h+m}(-n)^{h} \text{.}\end{align*}
	That is,
	\begin{align*} {n\brack m}&= (-1)^{n-m}\sum_{h=1}^{n-m}\binom{h+m-1}{m-1}{n\brack h+m}(-n)^{h} \text{.}\end {align*}
	To finish the proof, we just need that Equation \eqref{recst3} be true for $m=n$, which is obvious.
	\qed
	
%%4444444444444444444444444444444444444444444444444444444444444444444444444444444444444444444444444444444444444444444444444444444444444444444444444444444444444444444444444444444444%
\section {Extended congruences for the harmonic numbers $G^{(j+1)}_{p-1}$}
\begin{thm}
 Let  $k \ge 0$ be an integer and $p$ a prime number. We have 
\begin{align}\label{ge10ee}   G^{(k)}_{p-1}&=(-1)^k \sum_{j\ge0} (-1)^j\binom {j+k}{j}G^{(k+j)}_{p-1}p^j \text{.}\end{align}
In particular, when $k=0$, we have 
\begin{align}\label{ge10een}    \sum_{j\ge0} (-1)^jG^{(j+1)}_{p-1}p^j =0\text{.}\end{align}
	\begin{proof} Letting $n=p$ a prime number, and $m=k+1$ in Equation \eqref{recst3}, and dividing throughout by $(p-1)!$ provides the desired result.
	\end{proof}
\end{thm}

\indent Recall \cite{Gy19} that when $k\ge1$, the generalized Harmonic numbers $H^{(k)}_{p-1}$ admit the following p-adically converging expansion:
\begin{align}\label{e10e}   H^{(k)}_{p-1}&=(-1)^k \sum_{j\ge0} \binom {j+k-1}{j}H^{(k+j)}_{p-1}p^j\text{.}\end{align}
\noindent It is interesting to point out the similarity of Equation \eqref{e10e} and Equation \eqref{ge10ee}, but also some differences. Contrary to Equation \eqref{e10e}, the sum on the right-hand side of Equation \eqref{ge10ee} is finite: it is actually limited to $j=p-1-k$. We also notice that the sign alternates in Equation \eqref{ge10ee} and that there is a slight difference in the binomial coefficient.
\ \\

\indent It is also known \cite{Washington98} that, for odd prime $p$,
\begin{align}\label{bern}\sum_{j \ge 0}\binom{j+2k}{2k}B_j H^{(2k+j+1)}_{p-1}(-p)^j&=0\text{,}\end{align}
\noindent the convergence of the series being understood $p$-adically. More precisely, when $p\ge 5$, the following identity was shown in \cite{Gy19}:  
\begin{align}\sum_{j=0 }^{2n+1}\binom{j+2k}{2k}B_j H^{(2k+j+1)}_{p-1}(-p)^j& \equiv 0 \pmod {p^{2n+3}} \text{.} \end{align}

\noindent Now, we look for an equation similar to Equation \eqref{bern}, but involving the cycle Stirling numbers.  In the case where $k=0$, we have, for $p \ge 5$
\begin{equation}\label{ee10biss}\sum_{j=0 }^{2n+1}B_j H^{(j+1)}_{p-1}(-p)^j \equiv 0 \pmod {p^{2n+3}} \text{.}\end{equation}
For the first values of $n$, $n=0,1,2...$,  these congruences read
\begin{align*}
H_{p-1}+\frac{p}{2}H^{(2)}_{p-1} &\equiv 0 \pmod {p^3} \text{ , } \\
H_{p-1}+\frac{p}{2}H^{(2)}_{p-1} + \frac{p^2}{6}H^{(3)}_{p-1}&\equiv 0 \pmod {p^5}\text{ , }  \\
H_{p-1}+\frac{p}{2}H^{(2)}_{p-1} + \frac{p^2}{6}H^{(3)}_{p-1}- \frac{p^4}{30}H^{(5)}_{p-1}&\equiv 0 \pmod {p^7} \text{ ..., respectively. }
&\text{}
\end{align*}
A clue for our search of an equation analogous to Equation \eqref{bern} involving the Stirling cycle numbers is obtained by making use of Equation \eqref{conv} in order to recursively compute $H^{(j+1)}_{p-1}$ as function of the $G^{(i+1)}_{p-1}$, with $i\le j$. Then substituting $H^{(j+1)}_{p-1}$ in the above congruences and finally reducing modulo $p^{2n+3}$ as much as possible, by accounting for any previous congruence involving the $G^{(i+1)}_{p-1}$. In doing so, the following congruences are found, valid for $p\ge 5$:
\begin{align*}
 G_{p-1}-pG_{p-1}^{(2)} &\equiv 0 \pmod {p^3} \text{\ ,}\\
 G_{p-1}-pG_{p-1}^{(2)}+ \frac{p^2}{2}G^{(3)}_{p-1}&\equiv 0 \pmod {p^5}\text{\ ,}\\
 G_{p-1}-pG_{p-1}^{(2)}+ \frac{p^2}{2}G^{(3)}_{p-1}-\frac{p^4}{6}G^{(5)}_{p-1}&\equiv 0 \pmod {p^7}\text{\ ,}\\
 G_{p-1}-pG_{p-1}^{(2)}+ \frac{p^2}{2}G^{(3)}_{p-1}-\frac{p^4}{6}G^{(5)}_{p-1}+\frac{p^6}{6}G^{(7)}_{p-1}&\equiv 0 \pmod {p^9} \text{\ ...  \ etc.}
\end{align*}
The calculations become increasingly laborious as $n$ increases, but we are able to guess that
\begin{align} \label{conj1}
\sum_{j=0}^{2n+1}(j+1)B_jG^{(j+1)}_{p-1}p^j &\equiv 0 \pmod {p^{2n+3}} \text{.} \end{align}
In an even broader generalization, we anticipate the following theorem which will be proved in the next section.

\begin{thm} \label{5.2}
If $i \ge 1$ is an integer and $p$ a prime number, then we have 
 \begin{align} \label{conj1b}\sum_{j\ge0}B_j\binom{j+2i-1}{j}G^{(j+2i-1)}_{p-1}p^j =0\end{align}
or equivalently, 
 \begin{align} \label{conj2} \sum_{j\ge0}B_j\binom{j+2i-1}{j}{p \brack j+2i}p^j =0\text{.} \end{align}
\end{thm}
\noindent {\bf Remark.} Equation \eqref{conj1b} looks very much like Equation \eqref{bern} but the sums in Equation \eqref{conj1b} and Equation \eqref{conj2} are actually finite.

%5555555555555555555555555555555555555555555555555555555555555555555555555555555555555555555555555555555555555555555555555555555555555555555555555555555555555555555555555555555555555555555555555%%%%%
\section {The aerated triangular array  {{$\bf\it{{\mathcal A_{n,k}}}$}} }

\noindent For performing numerical verifications of Equation \eqref{conj2}, we now introduce the number $\mathcal A_{n,k}$.\\
\ \\
{\bf Definition.}  For non-negative integers $n$ and $k$, we define the number $\mathcal A_{n,k}$ by
\begin{equation} \label{def} \mathcal A_{n,k} := \sum_{h\ge0} B_h{k+h-1\choose h}{n\brack h+k}n^h . \end{equation}
It is clear from this definition that $\mathcal A_{n,k}$ is zero when $k>n$ and that  $\mathcal A_{n,n}=1$. The first terms of the sequence $\mathcal A_{n,k}$ are computed numerically and displayed in the following table.
{\setlength{\tabcolsep}{3pt}
{\center {\small
	\begin{tabular}{|c||ccccccccccc|}
	\hline
	$n$& $\mathcal A_{n,0}$& $\mathcal A_{n,1}$ & $\mathcal A_{n,2}$ & $\mathcal A_{n,3}$ & $\mathcal A_{n,4}$ & $\mathcal A_{n,5}$ & $\mathcal A_{n,6}$ & $\mathcal A_{n,7}$ & $\mathcal A_{n,8}$ & $\mathcal A_{n,9}$&$\mathcal A_{n,		10}$\\
	\hline
	\hline
	$0$  &$1$ & $0$ & $0$ &$0$& $0$ & $0$ & $0$&$0$&$0$&$0$&$0$\\
	$1$ &$0$ &$1$ & $0$ & $0$ &$0$& $0$ & $0$ & $0$&$0$&$0$&$0$\\
	$2$ &$0$&$0$ & $1$ & $0$ &$0$& $0$ & $0$ & $0$&$0$&$0$&$0$\\
	$3$&$0$&-1 & $0$ & $1$ &$0$& $0$ & $0$ & $0$&$0$&$0$&$0$\\
	$4$&$0$ &$0$ & -5 & $0$ &$1$& $0$ & $0$ & $0$&$0$&$0$&$0$\\
	$5$ &$0$&$24$ & $0$ & -15 &$0$& $1$ & $0$ & $0$&$0$&$0$&$0$\\
	$6$&$0$&$0$ & $238$ & $0$ &-35& $0$ & $1$ & $0$&$0$&$0$&$0$\\
	$7$ &$0$&-3396 & $0$ & $1281$ &$0$& -70 & $0$ & $1$&$0$&$0$&$0$\\
	$8$&$0$ &$0$ & -51508 & $0$ &$4977$& $0$ & -126 & $0$&$1$&$0$&$0$\\
	$9$ &$0$&$1706112$ & $0$ & -408700 &$0$& $15645$ & $0$ & -210&$0$&$1$&$0$\\
	$10$&$0$ &$0$ & 35028576&$ 0$&-2267320&$ 0$&$ 42273$&$0$&-330&$ 0$&$ 1$\\
	\hline
	\end{tabular}} 
\endcenter}
{\center {\bf Table 1:}  The triangular array $\mathcal A_{n,k}$, for $0 \le n,k \le 10$.\endcenter}

\noindent It is striking that these numbers seem to be zero when $n-k$ is odd, which, if true, would imply the validity of Theorem \ref{5.2}. This will be demonstrated in the next theorem. It is also striking that they seem to be all integers. 

\begin{thm}\label{n-kodd}
Let $n,k$ be non-negative integers. We have  $\mathcal A_{n,k}=(-1)^{n-k}\mathcal A_{n,k}$.
Equivalently, $\mathcal A_{n,k}=0$ when $n-k$ is odd.
	\begin{proof} When $n=0$, this is obviously true. Supposing $n>0$, we have
	\begin{align*}
	\mathcal A_{n,k}&= \sum_h B_h{k+h-1\choose h}{n\brack h+k}n^h\\
	&=(-1)^n \sum_h \frac{B_h}{n^k}{k+h-1\choose h}(-1)^{n-(h+k)}{n\brack h+k}(-n)^{h+k}\\
	&=(-1)^n \sum_h \frac{B_h}{n^k}{k+h-1\choose h}\sum_g {g-1\choose h+k-1}{n\brack g}(-n)^g    \text{ \ \ \  by Equation \eqref{recst4}}.\end{align*}
	But, it is easy to see that ${k+h-1\choose h}{g-1\choose h+k-1}={g-k\choose h}{g-1\choose k-1}$, so that
	\begin{align*}
	\mathcal A_{n,k}&= (-1)^{n-k}\sum_g  \sum_h B_h {g-k\choose h} {g-1\choose k-1}{n\brack g}(-n)^{g-k}  \\
	&=  (-1)^{n-k}\sum_g  (-1)^{g-k}B_{g-k}{g-1\choose k-1}{n\brack g}(-n)^{g-k} \text{\ \ \ \  by Equation \eqref{e2v}}\\
	&=  (-1)^{n-k}\sum_g  B_{g}{k+g-1\choose g}{n\brack g+k}n^{g} = (-1)^{n-k}\mathcal A_{n,k}. \end{align*} 
	\end{proof}
\end{thm}

\begin{thm} \label{int}
For non-negative integers $n$ and $k$, we have $ \mathcal A_{n,k} \in {\mathbb Z} $ .
\end{thm}

\noindent The proof of this theorem will not be given until the next section.  
Before proceeding to this proof, we want to point out a similarity between the triangle  $\mathcal  A_{n,k}$ and the triangle of Stirling numbers of the first kind.  We have the following theorem which is analogous to Lemma \ref{lemn-k}.

\begin{thm} 
\label{8.3}Let $n, k$ be non-negative integers such that $0\le k \le n$. There exists a polynomial  $P_k \in {\mathbb Q}[X]$, of degree $2k$, such that $P_k(n)$ coincides with $\mathcal A_{n,n-k}$. 
Moreover, when $k>0$, we have that $-1,0,\cdot \cdot \cdot,k$ are $k+2$ roots of $P_k(x)$.
	\begin{proof} By definition of $\mathcal{A}_{n,k}$, we have 
	\begin{equation*} \label{} \mathcal A_{n,n-k} = \sum_{h=0}^k B_h{n-1-(k-h)\choose h}{n\brack n-(k-h)}n^h,  \end{equation*}
	where the binomial coefficient is a polynomial in $n$ from ${\mathbb Q}[X]$, of degree $h$, and after Lemma \ref{lemn-k},  ${n\brack n-(k-h)}$ is a polynomial in $n$ from ${\mathbb Q}[X]$, of degree $2(k-h)$. 
	Therefore $\mathcal A_{n,n-k}$ is also a polynomial in $n$ from ${\mathbb Q}[X]$, of degree $2k$. Let $k>0$, and recall $Q_j $, the polynomial such that $Q_j(n)={n\brack n-j}$. We have 
	\begin{align*}
	P_k(u)&= \sum_{h=0}^k B_h{u+h-1-k\choose h}Q_{k-h}(u)u^h\\
	&= \sum_{h=0}^k B_h\frac{(u+h-1-k)\cdot\cdot(u-k)}{h!}Q_{k-h}(u)u^h.
  	\end{align*}
	From Lemma \ref{lemn-k}, we know that if $k>h$ then $0,1,\cdot \cdot \cdot,k-h$ are roots of the polynomial function $Q_{k-h}(x)$. Moreover when $k\ge h$, $Q_{k-h}(-1)=1$. \\
	\noindent Then, when $u=0$, we have $P_k(0)= B_0{-1-k\choose 0}Q_{k}(0)0^0=0 $, since $k>0$. When $u>0$, we have
	\begin{align*}
	P_k(u)&=\sum_{h=k-u+1}^{k} B_h\frac{(u-k)\cdot\cdot(u-k+h-1)}{h!}Q_{k-h}(u)u^h\\
	&=\sum_{h=k-u+1}^k (-1)^hB_h\frac{(k-u)\cdot\cdot(k-u-h+1)}{h!}Q_{k-h}(u)u^h.
 	\end{align*}
	 If $0<u\le k$, for any $h$ in the set $\{k-u+1,\cdot \cdot,k\}$ the product $(k-u)\cdot\cdot(k-u-h+1)$ must vanish because it has one factor which is zero, and then we also have $P_k(u)=0$.\\
	\noindent Finally, if $u=-1$,
	\begin{align*} \label{} 
	P_k(-1)&= \sum_{h=0}^k B_h\frac{(h-2-k)\cdot\cdot(-1-k)}{h!}(-1)^hQ_{k-h}(-1)\\
	&= \sum_{h=0}^k B_h\frac{(k+2-h)\cdot\cdot(k+1)}{h!}= \sum_{h=0}^k {k+1\choose h}B_h \\
	&= \sum_{h=0}^{k+1} {k+1\choose h}B_h - B_{k+1} =0 \text{ \ \  by Equation \eqref{e2v}.} 
	 \end{align*}
	\end{proof}
\end{thm}

%666666666666666666666666666666666666666666666666666666666666666666666666666666666666666666666666666666666666666666666666666666666666666666666666666666666666666666666666666666666666666666666666666666666%
\section {The dual triangle {{$\bf\it{{\mathcal B_{n,k}}}$}}}
\noindent We now introduce a dual triangle $\mathcal  B_{n,k}$ which is similar to the triangle of Stirling numbers of the second kind. Coming back to the polynomial $P_k$, we may extend the definition of $\mathcal A_{n,k}$ to non-positive indices since, for non-negative $n$ and $k$, it is natural to define $\mathcal A_{-n,-n-k}:= P_k(-n)$. Then, we have  
\begin{align*} \label{} \mathcal A_{-n,-n-k}&= \sum_{h=0}^k B_h{-n-1-(k-h)\choose h}Q_{k-h}(-n)(-n)^h\\
&= \sum_{h=0}^k B_h{n+k\choose h}Q_{k-h}(-n)n^h.\end{align*}
That is, 
\begin{equation*} 
\mathcal A_{-n,-k} =\sum_{h=0}^{k-n} B_h{k\choose h}Q_{k-n-h}(-n)n^h
\end{equation*}
or, accounting for Lemma \ref{lemn-k},
\begin{align*} 
\mathcal A_{-n,-k} &=\sum_{h=0}^{k-n} B_h{k\choose h}{k-h\brace n}n^h.\end{align*}
{\bf Definition.}  Let $n,k$ be positive integers. We define the number $\mathcal B_{n,k}$ by
\begin{equation} \label{def2} \mathcal B_{n,k} := \sum_{h\ge0} B_h{n\choose h}{n-h\brace k}k^h . \end{equation}
It is then clear that for all integers $n,k$, positive or negative, we have the duality:
\begin{equation} \label{8.5} \mathcal A_{-n,-k}= \mathcal B_{k,n} \text{.}   \end{equation}
This is similar to the duality $ {-n\brack -k}={k\brace n}   $ which holds \cite{Graham94} for the usual Stirling numbers. It is also clear from this definition that $\mathcal B_{n,k}$ is zero when $k>n$ and that  $\mathcal B_{n,n}=1$. Also note that we have $\mathcal B_{x+n,x}=P_n(-x)$. The first $\mathcal B_{n,k}$ are computed numerically and displayed in the following table.

{\setlength{\tabcolsep}{3pt}
{\center {\small
	\begin{tabular}{|c||cccccccccccc|}
	\hline
	$n$& $\mathcal B_{n,1}$ &$\mathcal B_{n,2}$ & $\mathcal B_{n,3}$ & $\mathcal B_{n,4}$ & $\mathcal B_{n,5}$ & $\mathcal B_{n,6}$ & $\mathcal B_{n,7}$ & $\mathcal B_{n,8}$ & $\mathcal B_{n,9}$ & $\mathcal B_{n,10}$&$\mathcal B_{n,11}		$&$\mathcal B_{n,12}$\\
	\hline
	\hline
	$1$  &$1$ & $0$  & $0$ & $0$ &$0$& $0$ & $0$ & $0$&$0$&$0$&$0$&$0$\\
	$2$  &$0$&$1$ & $0$ & $0$ &$0$& $0$ & $0$ & $0$&$0$&$0$&$0$&$0$\\
	$3$ &$0$ &$0$ & $1$ & $0$ &$0$& $0$ & $0$ & $0$&$0$&$0$&$0$&$0$\\
	$4$ &$0$ &-1 & $0$ & $1$ &$0$& $0$ & $0$ & $0$&$0$&$0$&$0$&$0$\\
	$5$ &$0$ &$0$ & -5 & $0$ &$1$& $0$ & $0$ & $0$&$0$&$0$&$0$&$0$\\
	$6$ &$0$ &$3$ & $0$ & -15 &$0$& $1$ & $0$ & $0$&$0$&$0$&$0$&$0$\\
	$7$&$0$&$0$ & $49$ & $0$ &-35& $0$ & $1$ & $0$&$0$&$0$&$0$&$0$\\
	$8$ &$0$&-17 & $0$ & $357$ &$0$& -70 & $0$ & $1$&$0$&$0$&$0$&$0$\\
	$9$ &$0$&$0$ & -809 & $0$ &$1701$& $0$ & -126 & $0$&$1$&$0$&$0$&$0$\\
	$10$ &$0$&$155$ & $0$ & -13175 &$0$& $6195$ & $0$ & -210 &$0$&$1$&$0$&$0$\\
	$11$ &$0$ &$0$ &$20317$&$ 0$&-120395&$ 0$&$ 18711$&$0$&-330&$ 0$&$ 1$&$0$\\
	$12$ & $0$ & -2073 & $0$ & $706893$ & $ 0$ &{\tiny -760100} & $ 0 $ & $ 49203$ & $0$ &  -495 & $0$ & $1$\\ 
	\hline 
	\end{tabular}} 
\endcenter}
{\center {\bf Table 2:}  The triangular array $\mathcal B_{n,k}$, for $1\le n,k \le 12$.\endcenter}

\noindent Again, we see on Table 2 that the $\mathcal B_{n,k}$  seem to be all integers and to vanish when $n-k$ is odd: this will be the next theorem. We can also, as was done in \cite{Graham94} for the
usual Stirling numbers, display  $\mathcal A_{n,k}$ and $\mathcal B_{-k,-n}$ in tandem, for $n,k \in \mathbb{Z}$. This is the purpose of Table 3, where we have left void the zero entries for $k>n$ and for odd $n-k$. The numbers which appear now in the diagonal lines are the values of the polynomial function $P_{n-k}(x)$  for integer arguments.

{\setlength{\tabcolsep}{3pt}
{\center {\small
	\begin{tabular}{|l||m{5mm}m{5mm}m{5mm}m{5mm}m{5mm}m{5mm}m{5mm}m{5mm}m{5mm}m{6mm}m{5mm}m{5mm}m{5mm}m{5mm}m{5mm}m{4mm}|}
	\hline
	 $ k \setminus n$ &-8 &-7 &-6 &-5 & -4 & -3 & -2 & -1 & 0 & 1 & 2 & 3&4&5&6&7\\
	\hline
	\hline
	-8  &1& &   &  &  & &  &  & &&&&&&&\\
	-7  &&$1$ &   &  &  & &  &  & &&&&&&&\\
	-6&-70 & &$1$ &  &  & &  &  & &&&&&&&\\
	-5 & &-35&&$1$ &  &  &&  &  & &&&&&& \\
	-4 &$357$&&-15 & & $1$&  &  &&  &  & &&&&&\\
	-3 &&$49$& &-5 &  & $1$ &  &&  &  & &&&&&\\
	-2 &-17&&$3$ & & -1 &  &$1$&  &&  &  & &&&&\\
	-1& &$0$& &$0$&  &$0$&& $1$ &  &&  &  & &&&\\
	0&$0$ &&$0$ & &$0$&  &$0$&  & $1$&  &&  &  & &&\\
	1& &$0$& &$0$ &&$0$  && $0$ &  & 1 &  &&  &  && \\
	2 &$0$ && $0$& &$0$  &  &$0$&  & $0$ & &$1$&  &&  &  & \\
	3&&$0$& &$0$ &  &$0$  &&  $0$ &  & -1&&$1$&  &&  &\\
	4&$0$&&$0$ & & $0$ &  &$ 0$& &$0$&&-5&&$ 1$&  &&  \\
	5&&$0$&& $0$& & $0$ &  & $0$ & & $ 24$ & &  -15 & & $1$&&\\ 
	6&$0$&& $0$& & $0$ &  &$0$&&$0$ &  & $238$ &  & -35 & &$1$& \\ 
	7&&0 & &0 &  & 0 & & 0 & & {\tiny-3396}   &  & 1281 &  & -70 & &1\\ 
	\hline 
	\end{tabular}} 
\endcenter}
{\center {\bf Table 3:}  $\mathcal A_{n,k}$ and $\mathcal B_{-k,-n}$ in tandem, for $-8 \le n,k  \le 7$.\endcenter}

\indent We now proceed to the proof of the integrality of these numbers.
	\begin{thm} \label{bnk}
	Let $n,k$ be non-negative integers, then $ \mathcal B_{n,k}$  is a triangular array of integers such that  $\mathcal B_{n,k}=0$ when $n-k$ is odd. Moreover, we have the inter-relations
	\begin{equation} \label{bnk1} \mathcal B_{x,x-n}= \sum_{u\ge 0}{n+x \choose n-u } {n-x \choose n+u} \mathcal A_{n+u,u}\end{equation}
	\begin{equation} \label{bnk2} \mathcal A_{x,x-n}= \sum_{u\ge 0}{n+x \choose n-u } {n-x \choose n+u} \mathcal B_{n+u,u}.\end{equation}
	\end{thm}

	\noindent {\bf Remark.} Theorem \ref{int} easily follows from Theorem \ref{bnk}. \\
	\noindent {\bf Remark.} Equation \eqref{bnk1} and Equation \eqref{bnk2} are formally the same as 
	\begin{align*}  {x \brace x-n}&= \sum_{u\ge 0}{n+x \choose n-u } {n-x \choose n+u} {u+n\brack u}\\
	{x \brack x-n}&=  \sum_{u\ge 0}{n+x \choose n-u } {n-x \choose n+u}{u+n\brace u} \text{,} \end{align*}
	 respectively, which hold \cite{Graham94} for the usual Stirling numbers.\\
 \ \\
	{\it{Proof of Theorem \ref{bnk}}.} For proving that $\mathcal B_{n,k}$ is integer, we are going to show that for any prime $p$, we have  $v_p(\mathcal B_{n,k})\ge 0$.\\ 
	\indent Firstly, we consider the case where $p$ divides $k$. For all $h$ such that 
	$h\ge 1$, we have $v_p\left(B_h{n\choose h}{n-h\brace k}k^h\right)\ge 0$, since $v_p(B_h) \ge -1$ by the Von Staudt-Clausen theorem and $v_p(k^h) \ge 1$. 
	Moreover, $v_p\left(B_0{n\choose 0}{n-0\brace k}k^0\right)=v_p\left({n\brace k}\right)\ge 0$, obviously. Then  $v_p(\mathcal B_{n,k})= v_p\left({n\brace k}+\sum_{h\ge1}B_h{n\choose h}{n-h\brace k}k^h \right)\ge 0$.\\
	\indent Then, we consider the case where $p$ does not divide $k$. We may write $$B_{n,k}=\sum_{h\ge0}B_h{n\choose h}{n-h\brace k}(k^h-1) + \sum_{h\ge0}B_h{n\choose h}{n-h\brace k}.$$ By the rule of multiplication of exponential generating 				functions \cite{Wilf92} and given the exponential generating functions Equation \eqref{egfs1} and Equation \eqref{egfbern},  we have 
	\begin{align*} \sum_{h\ge0}B_h{n\choose h}{n-h\brace k} &= [[x^n]] \left( \frac{x}{e^x -1} \frac{(e^x-1)^k}{k!} \right)\\&=\frac{1}{k} [[x^{n-1}]] \left(  \frac{(e^x-1)^{k-1}}{(k-1)!} \right)\\&=\frac{1}{k}{n-1\brace k-1}. \end{align*}
	But $v_p\left( \frac{1}{k}{n-1\brace k-1} \right) \ge 0$, since $p$ does not divide $k$, and then it suffices to show that for all values of $h$, $v_p\left(B_h(k^h-1) \right) \ge 0$. This is true because either $p-1$ does not divide $h$, and then by the 			Von Staudt-Clausen theorem, $v_p(B_h) \ge 0$, or $p-1$ divides $h$ and then $v_p(B_h) =-1$. But in this case, where $p$ does not divide $k$ and $p-1$ divides $h$, Fermat's little theorem holds and $v_p\left(k^h-1 \right) \ge 1$. \\
\ \\
	\indent Now, we turn to the proof of the inter-relations Equation \eqref{bnk1} and Equation \eqref{bnk2}. Let $p_n$ be a polynomial of degree $n$ from ${\mathbb Q}[X]$. The set of binomial coefficents $\left \{{x\choose k}; 0\le k \le n\right \}$ forms  a basis for the vector space of the polynomials from  ${\mathbb Q}[X]$ of degree less than $k+1$ and therefore there exists $a_{n,k}$ such that $p_n(x)=\sum_{k=0}^n a_{n,k}{x\choose k}$. 
	By the inversion formula Equation \eqref{hs1}, it is easy to verify that $a_{n,k}=\sum_u (-1)^{k-u}{k\choose u}p_n(u)$. We have seen that $\mathcal A_{x,x-n}= P_n(x)$ where $P_n$ is a polynomial of degree $2n$ from ${\mathbb Q}[X]$, so we can apply the above general inversion scheme to $P_n(-x)= \mathcal B_{x+n,x}$  and we obtain
	\begin{align*} \mathcal B_{x+n,x}&= \sum_{k=0}^{2n} \sum_{u=0}^{k}(-1)^{k-u} {k \choose u}\mathcal A_{u,u-n}  {-x \choose k}\\
	&= \sum_{k=n}^{2n} \sum_{u=n}^{k}(-1)^{u}{k \choose u} {x+k-1 \choose k}  \mathcal A_{u,u-n} \end{align*}
	
	\noindent which, given Theorem \ref{n-kodd}, clearly shows that $B_{n,k}=0$ when $n-k$ is odd. Now, since ${k \choose u} {x+k-1 \choose k}={x+k-1 \choose k-u} {x+u-1 \choose u}$, we have
	
	\begin{align*} \mathcal B_{x+n,x}&=\sum_{k=n}^{2n} \sum_{u=n}^{k}(-1)^{u}{x+k-1 \choose k-u} {x+u-1 \choose u} \mathcal A_{u,u-n}\\
	&= \sum_{u=n}^{2n}(-1)^{u}\sum_{k=u}^{2n}{x+k-1 \choose x +u-1} {x+u-1 \choose u} \mathcal A_{u,u-n}\\
	&= \sum_{u=n}^{2n}(-1)^{u}{2n+x \choose 2n-u} {x+u-1 \choose u} \mathcal A_{u,u-n}\\
	&= \sum_{u=n}^{2n}{2n+x \choose 2n-u} {-x \choose u} \mathcal A_{u,u-n}.\end{align*}
	
	Then, replacing $x$ by $x-n$ we have 
	\begin{align*} \mathcal B_{x,x-n}&= \sum_{u=n}^{2n}{n+x \choose 2n-u} {n-x \choose u} \mathcal A_{u,u-n}\\
	&= \sum_{u=0}^{n}{n+x \choose n-u} {n-x \choose n+ u} \mathcal A_{u+n,u}. \end{align*}
	Similarly, let $R_n(x)=\mathcal B_{x,x-n}$, which is a polynomial of degree $2n$ from ${\mathbb Q}[X]$. We apply the inversion to $R_n(-x)=\mathcal A_{x+n,x}$ which gives the similar identity where $\mathcal A$ and  $\mathcal B$ are 					exchanged, and this completes the proof of the theorem. 
	\qed 
\ 
\ \\

\indent For the usual Stirling numbers, there exist polynomials $\sigma_n(x)$ of degree $n-1$ from $\mathbb Q [X]$ (also known as Stirling polynomials in \cite{Graham94}, \cite{Knuth92}), such that 
\begin{equation}\label{kn1}  Q_n(x)={x \brack x-n}=x(x-1)\cdot\cdot\cdot(x-n)\sigma_n(x)   \end{equation}
or equivalently 
\begin{equation}\label{kn2}  Q_n(-x)={x+n \brace x}=(-1)^{n+1}x(x+1)\cdot\cdot\cdot(x+n)\sigma_n(-x).  \end{equation}
Similarly, accounting for Theorem \ref{8.3}, we can define the polynomial $\mathcal S_n(x)$ of degree $n-2$,  such that
\begin{equation}\label{kn3}  P_n(x)=\mathcal A_{x,x-n}=(x+1)x(x-1)\cdot\cdot\cdot(x-n)\mathcal S_n(x)   \end{equation}
or equivalently 
\begin{equation}\label{kn4}  P_n(-x)=\mathcal B_{x+n,x}=(-1)^{n+2}(x-1)x(x+1)\cdot\cdot\cdot(x+n)\mathcal S_n(-x).   \end{equation}
In the following table, we give the first instances of $\mathcal S_n(x)$, together with the Stirling polynomial $\sigma_n(x)$, given in  \cite{Graham94}.

{\setlength{\tabcolsep}{3pt}
{\center {\small
	\begin{tabular}{|c||cc c |cc c| cc c| c c|}
	\hline
	$n$ & & $1$ &  & & $2$& && $3$& && $4$\\
	\hline
	&&  & && & & &&&& \\
	$\sigma_n(x)$ && $\frac{1}{2}$ &&& $\frac{1}{24}(3x-1)$& && $\frac{1}{48}(x^2-x)$& & &$\frac{1}{5760}(15x^3-30x^2+5x+2)$ \\
	&&  & && & & &&&& \\
	$\mathcal S_{n}(x)$ && $0$&& & $-\frac{1}{24}$&& & $0$&& & $\frac{1}{5760}(7x^2+3x+2)$\\
	&&  & && & & &&&& \\
	\hline 
	\end{tabular}} 
\endcenter}
{\center {\bf Table 4:}   $\sigma_n(x)$ and $\mathcal S_n(x)$  for $n$ in the range $1$ to $4$.\endcenter}

\ \\
\indent The second and third columns of Table 2 are known to the OEIS \cite{OEIS}. Up to the sign and discarding the zeros, we find in these columns  the even index {\it Genocchi numbers} $G_{2n}$ and {\it Glaisher's G-numbers}, A001469 and A002111 at the OEIS, respectively, for which exponential generating functions are known. More generally, we have the following exponential generating function for  $\mathcal B_{n,k}$.

\begin{thm} \label{gnf}
Let $n,k$ be non-negative integers. We have
\begin{equation} \label{gnf1} \sum_{n \ge 0} \mathcal B_{n,k} \frac{x^n}{n!}=\frac {(e^x-1)^k}{k!}\frac{kx}{e^{kx}-1}.\end{equation}

	\begin{proof}
	The proof is straightforward: we use the rule of multiplication of exponential generating functions \cite{Wilf92}  and we have
	\begin{align*} \frac {(e^x-1)^k}{k!}\frac{kx}{e^{kx}-1}&= \left(\sum_{j\ge 0} {j \brace k}\frac{x^j}{j!} \right)\left(\sum_{j\ge 0} B_j\frac{(kx)^j}{j!}  \right)\\
	&=\sum_{u \ge 0}\left(\sum_{j+h=u}B_jk^j{u \choose h}{h \brace k}\right)\frac{x^u}{u!} \\
	&=\sum_{u \ge 0}\left(\sum_{j\ge 0}B_j{u \choose j}{u-j \brace k}k^j\right)\frac{x^u}{u!} . \end{align*}
	\end{proof}

\end{thm}

\noindent Unfortunately, the derivation of a similar generating function for $\mathcal A_{n,k}$ seems much more difficult.\\
\ \\
\indent We finish this section by pointing out a notable difference with the usual Stirling numbers. Whereas it is well-known that the usual Stirling matrices of both kinds are inverses of eachother, this is not the case for the $\mathcal A_{n,k}$  and $\mathcal B_{n,k}$. Let  $\mathcal A'_{n,k}$ be the inverse of the matrix $\mathcal A_{n,k}$. The first entries of  $\mathcal A'_{n,k}$ are displayed in the following table.

{\setlength{\tabcolsep}{3pt}
{\center 
	\begin{tabular}{|c||cccccccccc|}
	\hline
	$n$& $\mathcal A'_{n,1}$ & $\mathcal A'_{n,2}$ & $\mathcal A'_{n,3}$ & $\mathcal A'_{n,4}$ & $\mathcal A'_{n,5}$ & $\mathcal A'_{n,6}$ & $\mathcal A'_{n,7}$ & $\mathcal A'_{n,8}$ & $\mathcal A'_{n,9}$&$\mathcal A'_{n,10}$\\
	\hline
	\hline
	$1$  &$1$ & $0$ & $0$ &$0$& $0$ & $0$ & $0$&$0$&$0$&$0$\\
	$2$ &$0$ & $1$ & $0$ &$0$& $0$ & $0$ & $0$&$0$&$0$&$0$\\
	$3$ &$1$ & $0$ & $1$ &$0$& $0$ & $0$ & $0$&$0$&$0$&$0$\\
	$4$ &$0$ & $5$ & $0$ &$1$& $0$ & $0$ & $0$&$0$&$0$&$0$\\
	$5$ &$-9$ & $0$ & $15$ &$0$& $1$ & $0$ & $0$&$0$&$0$&$0$\\
	$6$&$0$ & $-63$ & $0$ &$35$& $0$ & $1$ & $0$&$0$&$0$&$0$\\
	$7$ &$1485$ & $0$ & $-231$ &$0$& $70$ & $0$ & $1$&$0$&$0$&$0$\\
	$8$ &$0$ & $18685$ & $0$ &$-567$& $0$ & $126$ & $0$&$1$&$0$&$0$\\
	$9$ &$-844757$ & $0$ & $125515$ &$0$& $-945$ & $0$ & $210$&$0$&$1$&$0$\\
	$10$ &$0$ &$-14862727$&$ 0$&$600655$&$ 0$&$-693$&$0$&$330$&$ 0$&$ 1$\\
	\hline
	\end{tabular}
  \endcenter}
  {\center {\bf Table 5:}  The triangular array $\mathcal A'_{n,k}$, for $1 \le n,k \le 10$.\endcenter}

\noindent  We do not see any evident link between the matrix $\mathcal A'_{n,k}$ and the matrix $\mathcal B_{n,k}$.
\noindent Moreover, whereas the Stirling matrices are convolution matrices in the sense of \cite{Knuth92}, this is not the case for the matrices  $\mathcal A_{n,k}$ and  $\mathcal B_{n,k}$. This is easily checked on their first entries, as indicated in \cite{Knuth92}.  If $\mathcal A_{n,k}$ and $\mathcal B_{n,k}$  were convolution matrices, inverse of eachother, there would exist a function $f(x)$ such that $\sum_n \mathcal B_{n,k}\frac{x^n}{n!}=\frac {f(x)^k}{k!}$ and the generating function for $\mathcal A_{n,k}$ would read $\sum_n \mathcal A_{n,k}\frac{x^n}{n!}=\frac {g(x)^k}{k!}$, where $g$ is the compositional inverse of $f$. The triangular array $\frac{n!}{k!}\mathcal B_{n,k}$ is not even a Riordan array for which (see for instance \cite{Luzon12}) the ordinary generating function reads $d(x)\cdot h(x)^k$ for some power series $d(x)$ and $h(x)$, with $h(0)=0$ and $Dh(0) \neq 0$. Knowing the generating function Equation \eqref{gnf1} for $\mathcal B_{n,k}$, the difficulty to find an analogous generating function  for $\mathcal A_{n,k}$ has to do with these observations.

%777777777777777777777777777777777777777777777777777777777777777777777777777777777777777777777777777777777777777777777777777777777777777777777777777777777777777777777777777777777777777777777777777%
\section {Discussion and Questions}

\indent  Since the second column of the triangle $\mathcal B_{n,k}$ corresponds to the Genocchi numbers, we might consider the other columns as some sort of generalized Genocchi numbers. However, these numbers are not the same as the already known generalized Genocchi numbers from \cite{Domaratzki04}, nor as those from \cite{Luo11}. For the classical Genocchi numbers, there exists a recursion so that  for $n \ge 1$, it is possible to compute  $\mathcal B_{2n,2}$ recursively: $$ \mathcal B_{2n,2} = n  - \frac{1}{2}\sum_{j=1}^{n-1} {2n \choose  2j}\mathcal B_{2j,2}.$$
Moreover, a combinatorial interpretation has been given to the Genocchi numbers \cite{Dumont74}, but to the author's knowledge, this is not the case for Glaisher's G-numbers, $ \left| \mathcal B_{2n+1,3}\right|$. Here we raise the more general questions: for a given  $k>2$,  find a recursion for the {\it Bernoulli-Stirling numbers of the second kind}  $\mathcal B_{n,k}$ and find combinatorial objects that they enumerate.\\

\indent As for the {\it Bernoulli-Stirling numbers of the first kind} $\mathcal A_{n,k}$,  we have even more questions. Apart from their appearance in the above investigation of congruences modulo prime powers for the cycle Stirling numbers, we don't know their mathematical interest. Any recurrence that would allow to compute an entry in this triangular array from entries of previous lines would be insightful, and might lead to a direct proof of Theorem \ref{int}. Moreover $\mathcal A_{n,k}$ cries for a generating function, of any kind, or at least a functional equation involving such a function. There also remains the problem of the combinatorial interpretation of $\mathcal A_{n,k}$.\\
\ \\
\indent By comparison to these quite complicated combinatorics questions, the study of the arthmetic properties of the Bernoulli-Stirling numbers would seem more easy, as it could be made use of their explicit expression in terms of Bernoulli and Stirling numbers and then take advantage the existing knowledge on the arithmetic properties of the latter. In particular, we might expect that the Bernoulli-Stirling numbers satisfy some sort of Kummer congruence, as do the regular Bernoulli numbers.

\end{document}